\documentclass[a4paper,12pt,reqno]{amsart}

\usepackage{amssymb,amsmath,amsthm}
\usepackage{esint,hyperref}
\usepackage{tikz}
\usetikzlibrary{decorations.pathreplacing}
\hypersetup{hidelinks}

\numberwithin{equation}{section}
\allowdisplaybreaks
\newcommand{\dd}{\mathop{}\!\mathrm{d}}

\newcommand{\Ric}{\operatorname{Ric}}

\newcommand{\II}{\operatorname{II}}

\newtheorem{thm}{Theorem}[section]

\newtheorem{prop}[thm]{Proposition}

\newtheorem{rem}[thm]{Remank}

\theoremstyle{definition}
\newtheorem{defn}[thm]{Definition}
\newtheorem*{question}{Question}

\title[Relative isoperimetric comparison]{Isoperimetric comparison for Surfaces with nonnegative curvature and uniformly convex boundary}

\author{Chao Xia}
\address{(C.X.) School of Mathematical Sciences,
Xiamen University, 361005, Xiamen, P.~R.~China}
\email{chaoxia@xmu.edu.cn}
\thanks{C.X. is supported by NSFC (Grant No. 12271449, 12671069, 12526203, 12526102) and the Natural Science Foundation of Fujian Province of China (Grant No. 2024J011008).}

\author{Jincong Zheng}
\address{(J.Z.) School of Mathematical Sciences,
Xiamen University, 361005, Xiamen, P.~R.~China}
\email{zhengjincong@stu.xmu.edu.cn}
\date{}

\begin{document}

\begin{abstract}
We establish a sharp relative isoperimetric comparison for compact surfaces with nonnegative Gaussian curvature and
smooth nonempty boundary whose geodesic curvature is at least $1$.
We prove that the normalized relative isoperimetric profile for such a surface is bounded
below by that of the Euclidean unit disk, with equality holding at some area
fraction in $(0,1)$ if and only if the surface is isometric to the unit disk.
The proof uses a differential comparison argument for the isoperimetric
profile, based on a suitable non-constant normal variation of isoperimetric
curves inspired by the model geometry.
\end{abstract}
\maketitle
\section{Introduction}
The celebrated Levy-Gromov isoperimetric inequality \cite{Gr} provides a sharp comparison between the isoperimetric profile of closed Riemannian $n-$manifolds with $Ric\ge n-1$ and that of unit sphere $\mathbb{S}^n$. 
The main ingredient of Gromov's proof is the Heintze-Karcher comparison \cite{HK} on the isoperimetric hypersurfaces. An alternative proof by Bayle \cite{Ba} is based on a differential inequality on the isoperimetric profile.

For compact Riemannian manifolds with boundary, 
the relative isoperimetric problem concerns minimizers for the perimeter of the relative
boundary of domains enclosing a prescribed
volume. This problem is also called the partitioning
problem. Bayle-Rosales \cite{BayleRosales} proves a Levy-Gromov-type sharp isoprimetric comparison between the relative isoperimetric profile of compact Riemannian manifolds with $Ric\ge n-1$ and convex boundary,  and that of unit half sphere $\mathbb{S}^n_+$,
by following the method of Bayle \cite{Ba}.

 The family of complete Riemannian manifolds with $Ric\ge 0$ and uniformly convex boundary is of great interest. 
Recently, Yan-Zhu \cite{YanZhu} proves that such a manifold with boundary must be compact, confirming a conjecture by M. Li \cite{Li}. See also Ge-Li \cite{GeLi} for the case of nonnegative sectional curvature.  Yan-Zhu's theorem can be viewed as a boundary version of Bonnet-Myers compactness theorem for closed manifolds with $Ric\ge n-1$. 
X. Wang \cite{Wang} made the conjecture that for a compact Riemannian manifold $M^n$ with $Ric\ge 0$ and the boundary second fundamental form $II_{\partial M}\ge 1$, its boundary area $|\partial M|\le |\mathbb{S}^{n-1}|$ with equality holding iff $M$ is isometric to a Euclidean unit ball $\mathbb{B}^n$. He confirmed it for $n=3$.

We will study relative isoperimetric comparison for compact Riemannian manifolds with $Ric\ge 0$ and uniformly convex boundary. We recall the definition of relative isoperimetric profile. Let $M$ be a smooth, compact, connected Riemannian manifold with smooth boundary $\partial M$. Let $|\cdot|$ denote Riemannian measure in the dimension of the set under consideration.
\begin{defn}[Relative isoperimetric profile]\label{def:profile}
 For $0\le a\le|M|$, set
\begin{equation}\label{eq:profile}
 I_M(a)=\inf\left\{
 |\partial^*E\cap\operatorname{int}M|:
 E\subset M\text{ has finite perimeter},\ |E|=a
 \right\},
\end{equation}
where $\partial^*E$ is the reduced boundary of $E$. The function
\begin{equation}\label{eq:normalized-profile}
 h_M(t)=\frac{I_M(t|M|)}{|M|},\qquad 0<t<1,
 \qquad h_M(0)=h_M(1)=0,
\end{equation}
is called the relative isoperimetric profile of $M$.
\end{defn}

The isoperimetric region of $\mathbb{B}^n$ is known to be the intersection of another ball with $\mathbb{B}^n$, which is due to Burago-Maz’ya \cite{BM}. Hence, the relative isoperimetric profile of $\mathbb{B}^n$ can be explicitly computed.
Inspired by the Levy-Gromov isoperimetric inequality for closed manifolds, we ask the following question: 
\begin{question}
Let $(M^n,g)$, $n\ge2$, be a compact connected Riemannian
manifold with smooth boundary. Suppose that
\[
 \Ric_M\ge0,\qquad \II_{\partial M}\ge g|_{T\partial M},
\]
Is it true that
\[
 h_M(t)\ge h_{\mathbb{B}^n}(t),\qquad 0\le t\le1?
\]
\end{question}



In this paper, we confirm the $2$-dimensional case affirmatively. The higher dimensional case for the above Question will be studied in a forthcoming paper.  For $(M^2 ,g)$, we use $K_M$ and $\kappa_{\partial M}$ to denote the Gaussian curvature of $M$ and the geodesic curvature of $\partial M$ respectively. Our main result is the following. 
\begin{thm}\label{thm:main}
Let $(M^2,g)$ be a compact, connected  Riemannian surface with smooth nonempty boundary. Suppose that
\[
K_M\ge0,\qquad \kappa_{\partial M}\ge 1.
\]
 Then
\[
 h_M(t)\ge h_{\mathbb{B}^2}(t),\qquad 0\le t\le1.
\]
Moreover, equality holds for some $t\in(0,1)$ if and only if
$M$ is isometric to $\mathbb{B}^2$.
\end{thm}

\medskip
\noindent\textbf{Proof strategy.} 
Our proof uses a differential comparison argument for the isoperimetric
profile. Variations of isoperimetric regions provide smooth upper
supports, reducing the profile comparison to estimates for the
quadratic form $Q$. In Euclidean space and Euclidean half-space, the classical
isoperimetric interfaces are spheres and hemispheres, respectively.
Constant normal variations preserve these families of isoperimetric
interfaces and therefore recover the exact model second variation.
In contrast, the isoperimetric interfaces in $\mathbb{B}^n$ are
free-boundary spherical caps, with the equatorial disk as a limiting
case, and the normal speed of this family is not constant along the
interface. This motivates using the model normal speed as a guide
for constructing suitable test functions. So we first compute the nonconstant normal speed of
the isoperimetric regions in the unit disk as the enclosed area varies.
This motivates corresponding normalized test functions on
free-boundary arcs of geodesic curvature $k$ and length $L$.
The curvature assumptions yield an explicit upper bound
$\widehat Q(k,L)$, which recovers the exact model second variation
when $L=L_0(k)$.

The main difficulty is that the curvature and length of a general
isoperimetric interface need not satisfy the model relation.
At a minimum of the profile difference, the first variation condition
matches the curvatures. The area bound and the minimum condition imply $L\le L_0(k)$.
The strict monotonicity of $L\widehat Q(k,L)$ in $L$ then allows
us to compare the second variation bound with its model value. In the proof of the main
theorem, the Gauss--Bonnet theorem and the Heintze--Karcher inequality
exclude a closed interface at this minimum, ensuring that the
free-boundary arc estimates apply. The second variation comparison
gives the desired inequality, while strict monotonicity and the
equality case of the Heintze--Karcher inequality yield rigidity.

This paper is organised as follows.  In Section~2, we collect some known properties for $M$ and its isoperimetric regions.  In Section~3, we study some properities of the isoperimetric profile of the model space
$\mathbb{B}^2$. In section~4, we use a suitable test function inspired by the model case to prove a sharp second variation comparison.  In Section~5, we prove the main theorem. 

\

\medskip
\noindent\textbf{Disclosure on AI assistance.}
A publicly available large language model (LLM) was used to suggest
the test function $u_{k,L}$, expedite the computations of standard quantities
in the model space, and assist in the verification of the correctness of the calculations. The manuscript was written entirely by the authors. The authors take full responsibility for the paper.

\

\section{Preliminaries}\label{sec:preliminaries}

We collect some known properties for $M$ and its isoperimetric regions, which are known to experts and will be used in this paper.
\begin{prop}\label{lem:geometry-two}
Let $M^2$ be a complete, connected surface with 
boundary $\partial M$. Assume $K_M\ge0$ and $\kappa_{\partial M}\ge1$. Then $M$ is compact and
diffeomorphic to the closed disk, and
\[
 \int_M K_M\,\dd s<2\pi,
 \qquad |M|\le\frac{|\partial M|}{2}\le\pi.
\]
There is no simple closed geodesic in $\operatorname{int}M$.
\end{prop}

\begin{proof} 
Compactness of $M$ is proved by 
Ge-Li \cite{GeLi}, see also \cite[Theorem~1.3]{Wu}.

By the Gauss--Bonnet theorem,
\[
 2\pi\chi(M)
 =\int_M K\,\dd s+\int_{\partial M}\kappa\,\dd s>0.
\]
Hence $M$ is diffeomorphic to the closed disk and $\chi(M)=1$. Thus
\[
 \int_M K\,\dd s
 =2\pi-\int_{\partial M}\kappa\,\dd s<2\pi.
\]
Moreover, since $K_M\ge0$ and $\kappa_{\partial M}\ge1$,
\[
 |\partial M|
 \le\int_{\partial M}\kappa\,\dd s
 =2\pi-\int_M K\,\dd s
 \le2\pi.
\]
The Heintze-Karcher inequality
\cite[Theorem~1]{Ros} gives
\[
 |M|\le\frac12\int_{\partial M}\frac1\kappa\,\dd s
 \le\frac{|\partial M|}{2}\le\pi.
\]

Finally, suppose that $\gamma$ is a simple closed geodesic in
$\operatorname{int}M$. Since $M$ is a disk, $\gamma$ bounds a disk
$D\subset\operatorname{int}M$. Applying the Gauss--Bonnet theorem to
$D$, we obtain
\[
 2\pi=\int_D K\,\dd s
 \le\int_M K\,\dd s<2\pi,
\]
a contradiction.
\end{proof}

\begin{prop}\label{prop:regularity}
Let $M$ be a compact, connected Riemannian surface with smooth boundary.
For every $0<a<|M|$, the infimum in
\eqref{eq:profile} is attained. If $E$ is an isoperimetric region with $|E|=a$, then
\[
 \gamma=\overline{\partial^*E\cap\operatorname{int}M}
\]
is smooth, compact and embedded, and
$|\gamma|=I_M(a)$. Moreover, $\gamma$ has constant geodesic
curvature with respect to the outward normal of $E$.
If $\gamma$ intersects $\partial M$, then it must be orthogonally.
Under the hypotheses of Proposition~\ref{lem:geometry-two},
$\gamma$ is connected and is either a simple closed curve
contained in $\operatorname{int}M$ or an arc with two endpoints
on $\partial M$.
The profile $h_M$ is continuous on $[0,1]$, positive on $(0,1)$, and
\[
 h_M(t)=h_M(1-t).
\]
\end{prop}

\begin{proof}
See \cite[Proposition~2.3 and Remarks~2.6 and~3.8]{BayleRosales}
for existence, regularity and the basic properties of the profile.

Under the hypotheses of Proposition~\ref{lem:geometry-two},
suppose that $\gamma$ is not connected.
By \cite[Proposition~3.9(ii) and Remark~2.6]{BayleRosales},
the assumptions $K_M\ge0$ and $\kappa_{\partial M}\ge1$
imply that $\gamma$ does not meet $\partial M$ and that
each of its connected components is a simple closed geodesic.
This contradicts Proposition~\ref{lem:geometry-two}.
Hence $\gamma$ is connected.
Being a smooth, compact, connected embedded curve,
it is either a simple closed curve contained in
$\operatorname{int}M$ or an arc with two endpoints on $\partial M$.
\end{proof}

Let $\nu$ and $N$ be the outward unit normal to $\partial M$ and $\gamma$, respectively. let $\eta$ be the outward unit co-normal to $\partial\gamma$ in $\gamma$. Let $\kappa_{\partial M}$ and $\kappa_\gamma$ be the geodesic curvature of $\partial M$ (with respect to $\nu$) and $\gamma$  (with respect to $N$) respectively.
For a free-boundary curve $\gamma$, $\eta=\nu$ on $\partial\gamma$.

\begin{prop}[First and second variations]\label{prop:variation}
Let $E$ be an isoperimetric region with its relative boundary $\gamma$ having constant
geodesic curvature $\kappa_\gamma$ and meets $\partial M$ orthogonally at $\partial \gamma$. For any 
$u\in C^\infty(\gamma)$, there is a smooth vector field $X$ on $M$,
tangent to $\partial M$, such that $X=uN$ on $\gamma$. Let $\varphi_\tau$
be its local flow, and set $E_\tau=\varphi_\tau(E)$ and
$\gamma_\tau=\varphi_\tau(\gamma)$. Then
\[
 \left.\frac{\dd}{\dd\tau}|E_\tau|\right|_{\tau=0}
 =\int_\gamma u\,\dd s,
 \qquad
 \left.\frac{\dd}{\dd\tau}|\gamma_\tau|\right|_{\tau=0}
 =\kappa_\gamma\int_\gamma u\,\dd s.
\]
 Normalize $u$ such that $\int_\gamma u \dd s=1$.  
If we view the length functional $L$ of $\gamma$ as a function of the area function $A$ of $E$, then
\[\frac{d L}{dA}\bigg|_{A=|E|}=\kappa_\gamma,\]
and
\begin{equation}\label{eq:Q}
 \frac{d^2 L}{dA^2}\bigg|_{A=|E|}=\int_\gamma\left(
 |u'(s)|^2
 -\bigl(\kappa_\gamma^2+K_M\bigr)u^2
 \right)\dd s
 -\sum_{p\in \partial \gamma} \kappa_{\partial M}(p)u(p)^2:=Q(u).
\end{equation}

\end{prop}

\begin{proof}
See \cite[Chapter~9, Theorems~9.7--9.8]{Ritore}.
\end{proof}

Note the isoperimetric profile $h_M$ is continuous, but might not be differentiable, it is helpful to have a smooth upper support function replacing $h_M$.
\begin{prop}\label{lem:upper-support}
Let $E$ be an isoperimetric region with
$|E|=t_0|M|$, where $0<t_0<1$, and  $\gamma$ be its relative
boundary with constant geodesic curvature $\kappa_\gamma$. For every $u\in C^\infty(\gamma)$ with
$\int_\gamma u\,\dd s=1$, there is a smooth function $\widehat h_M$ near
$t_0$ such that
\begin{equation}\label{eq:upper-support}
 \begin{gathered}
 \widehat h_M(t)\ge h_M(t),\qquad \widehat h_M(t_0)=h_M(t_0),\\
 \widehat h_M'(t_0)=\kappa_\gamma,
 \qquad \widehat h_M''(t_0)=|M|Q(u).
 \end{gathered}
\end{equation}
\end{prop}

\begin{proof}
Apply the construction in
\cite[proof of Theorem~3.2, equations~(3.5)--(3.6)]{BayleRosales}
to the variation with normal speed $u$, and use Proposition \ref{prop:variation}.
\end{proof}





\section{The isoperimetric profile of the model space
\texorpdfstring{$\mathbb{B}^2$}{B2}}

For $k>0$, let $E_k\subset \mathbb{B}^2$ be the smaller region bounded by
$\partial \mathbb{B}^2$ and a circular arc of radius $1/k$ meeting $\partial \mathbb{B}^2$
orthogonally. Let $E_0$ be a half-disk. They are isoperimetric regions of $\mathbb{B}^2$. Let $\gamma_k$ be the relative
boundary of $E_k$.   The geodesic curvature of $\gamma_k$ is $\kappa_\gamma=k$.
By \cite[Theorem~9.26]{Ritore}, the regions $E_k$, $k\ge0$, and their
complements are precisely the relative isoperimetric regions in $\mathbb{B}^2$,
up to rotations.

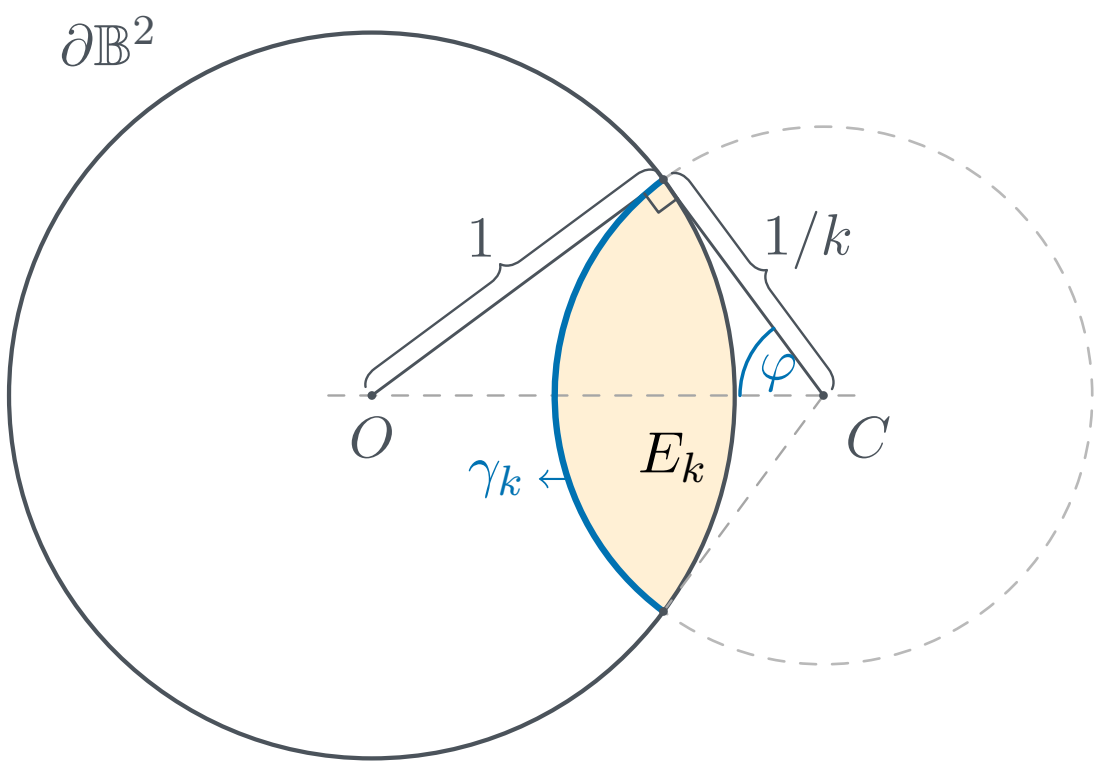
\begin{figure}[htbp]
\centering
\begin{tikzpicture}[scale=2.4, font=\small, line cap=round, line join=round]
 \definecolor{modelblue}{HTML}{0072B2}
 \definecolor{modelgray}{HTML}{4B535B}
 \definecolor{modelshade}{HTML}{FFF0D5}
 \pgfmathsetmacro{\capradius}{1/1.35}
 \pgfmathsetmacro{\centerd}{sqrt(1+\capradius*\capradius)}
 \pgfmathsetmacro{\capangle}{atan(1.35)}
 \pgfmathsetmacro{\capstart}{180-\capangle}
 \pgfmathsetmacro{\capend}{180+\capangle}
 \pgfmathsetmacro{\px}{1/\centerd}
 \pgfmathsetmacro{\py}{\capradius/\centerd}
 \def\marksize{0.065}
 \coordinate (O) at (0,0);
 \coordinate (C) at (\centerd,0);
 \coordinate (P) at (\px,\py);
 \coordinate (Pm) at (\px,-\py);
 \begin{scope}
  \clip (O) circle[radius=1];
  \fill[modelshade] (C) circle[radius=\capradius];
 \end{scope}
 \draw[gray!55,dashed,line width=0.5pt] (C) circle[radius=\capradius];
 \draw[modelgray,line width=0.8pt] (O) circle[radius=1];
 \draw[gray!70,dashed,line width=0.45pt] (-0.12,0)--({\centerd+0.12},0);
 \draw[modelgray,line width=0.55pt] (O)--(P);
 \draw[modelgray,line width=0.55pt] (C)--(P);
 \draw[modelgray,line width=0.45pt,decorate,
  decoration={brace,amplitude=3.5pt,raise=2pt}]
  (O)--(P);
 \node[modelgray] at
  ({\px/2-0.17*\py},{\py/2+0.17*\px}) {$1$};
 \draw[modelgray,line width=0.45pt,decorate,
  decoration={brace,amplitude=3.5pt,raise=2pt}]
  (P)--(C);
 \node[modelgray] at
  ({(\px+\centerd)/2+0.22*\py/\capradius},
   {\py/2+0.22*(\centerd-\px)/\capradius}) {$1/k$};
 \draw[gray!70,dashed,line width=0.45pt] (C)--(Pm);
 \draw[modelblue,line width=1.2pt]
  (P) arc[start angle=\capstart,end angle=\capend,radius=\capradius];
 \draw[modelgray,line width=0.5pt]
  ({\px-\marksize*\px},{\py-\marksize*\py})
  --({\px-\marksize*\px+\marksize*(\centerd-\px)/\capradius},
     {\py-\marksize*\py-\marksize*\py/\capradius})
  --({\px+\marksize*(\centerd-\px)/\capradius},
     {\py-\marksize*\py/\capradius});
 \draw[modelblue,line width=0.65pt]
  ({\centerd+0.23*cos(\capstart)},{0.23*sin(\capstart)})
  arc[start angle=\capstart,end angle=180,radius=0.23];
 \node[modelblue] at
  ({\centerd-0.14*cos(\capangle/2)},{0.14*sin(\capangle/2)}) {$\varphi$};
 \node[modelblue,anchor=east,inner sep=1pt] at (0.44,-0.23) {$\gamma_k$};
 \draw[modelblue,line width=0.4pt,->]
  ({\centerd-sqrt(\capradius*\capradius-0.23*0.23)},-0.23)--(0.46,-0.23);
 \node at (0.83,-0.17) {$E_k$};
 \node[modelgray,above left] at (-0.53,0.85) {$\partial \mathbb{B}^2$};
 \foreach \point in {O,C,P,Pm}
  \fill[modelgray] (\point) circle[radius=0.012];
 \node[below,modelgray] at (O) {$O$};
 \node[below right,modelgray] at (C) {$C$};
\end{tikzpicture}
\caption{The model region $E_k$ for $k>0$, with boundary $\gamma_k$ in blue.
The centers of the unit and auxiliary circles are $O$ and $C$,
respectively. The angle at $C$ is $\varphi$.}
\label{fig:model-regions}
\end{figure}

For $k>0$, let $\varphi\in(0,\pi/2)$ be half the central angle of the
circular arc in Figure~\ref{fig:model-regions}, measured at the center
$C$ of the auxiliary circle. Orthogonality gives $k=\tan\varphi$.
Set
\[
 A_0(k)=|E_k|,\qquad L_0(k)=|\gamma_k|.
\]

\begin{prop}\label{prop:disk-model}
The relative isoperimetric regions in $\mathbb{B}^2$ are the regions $E_k$ and
their complements, up to rotations. For $k>0$,
\begin{equation}\label{eq:disk-volume-area}
 \begin{aligned}
 A_0(k)&=\frac\pi2-\arctan k+\frac{\arctan k}{k^2}-\frac1k,\\
 L_0(k)&=\frac{2\arctan k}{k}.
 \end{aligned}
\end{equation}
Both functions are strictly decreasing for $k>0$, and extend smoothly
to $k=0$, with
\[
 A_0(0)=\frac\pi2,\qquad L_0(0)=2.
\]
The isoperimetric profile $ h_{\mathbb{B}^2}(t), t\in (0,\frac12]$ of $\mathbb{B}^2$ satisfies
\begin{equation}\label{eq:disk-profile-derivatives}
 \begin{aligned}
 h_{\mathbb{B}^2}(t)&=\frac{L_0(k)}\pi,\qquad h_{\mathbb{B}^2}'(t)=k,\\
 h_{\mathbb{B}^2}''(t)&=-\frac{\pi k^2}{L_0(k)-2/(1+k^2)}\quad t\in (0,\frac12),\\
 h_{\mathbb{B}^2}''(1/2)&=-\frac{3\pi}{4}\quad t=\frac12,
 \end{aligned}
\end{equation}
where $k=k(t)$ is such that $A_0(k(t))=|\mathbb{B}^2| t=\pi t$. 
In particular, $h_{\mathbb{B}^2}\in C^2(0,1)$.
\end{prop}

\begin{proof}
The computation of $A_0$ and $L_0$ is direct, see for example
\cite[the proof of Proposition~1]{AndrewsBryan}. Rewriting their
expressions in terms of $k$ gives \eqref{eq:disk-volume-area}.
A simple calculation gives
\[
 A_0'(k)
 =\frac{2}{k^2(1+k^2)}-\frac{2\arctan k}{k^3}
 =-\frac{L_0(k)-2/(1+k^2)}{k^2},
 \qquad L_0'(k)=kA_0'(k).
\]
Since
\[
 \arctan k-\frac{k}{1+k^2}
 =\int_0^k\frac{2x^2}{(1+x^2)^2}\,\dd x>0
 \qquad(k>0),
\]
it follows that $A_0'(k)<0$, and hence both $A_0$ and $L_0$ are strictly
decreasing.

The Taylor series of $\arctan k$ shows that both functions extend smoothly
to $k=0$. More explicitly,
\[
 A_0(k)=\frac\pi2-\frac43k+\frac{8}{15}k^3+O(k^5),
 \qquad L_0(k)=2-\frac23k^2+O(k^4).
\]
Thus $A_0(0)=\pi/2$, $L_0(0)=2$, and $A_0'(0)=-4/3$.

Since $A_0$ is strictly decreasing from $\pi/2$ to $0$, each
$t\in(0,1/2)$ determines a unique $k=k(t)>0$ such that $A_0(k)=\pi t$. Thus
\[
 h_{\mathbb{B}^2}(t)=\frac{L_0(k)}\pi.
\]
Note $\frac{dk}{dt}=\frac{\pi}{A_0'(k)}$. Then
the chain rule gives
\[
 h_{\mathbb{B}^2}'(t)=\frac{L_0'(k)}{A_0'(k)}=k,
\]
\[
 h_{\mathbb{B}^2}''(t)=\frac{dk}{dt}=\frac\pi{A_0'(k)}
 =-\frac{\pi k^2}{L_0(k)-2/(1+k^2)}.
\]
As $t\uparrow1/2$, we have $k\downarrow0$, so
\[
 h_{\mathbb{B}^2}'(t)\longrightarrow0,
 \qquad h_{\mathbb{B}^2}''(t)\longrightarrow-\frac{3\pi}{4}.
\]
Together with the symmetry $h_{\mathbb{B}^2}(t)=h_{\mathbb{B}^2}(1-t)$, these limits show
that $h_{\mathbb{B}^2}\in C^2(0,1)$ and $h_{\mathbb{B}^2}''(1/2)=-3\pi/4$.
\end{proof}

     For $k\ge 0$, parametrize the model arc $\gamma_k$ by its arclength $s$, i.e.
\[
 \tilde{X}(k, s)=
 \left(\frac{\sqrt{1+k^2}-\cos ks}{k},
 \frac{\sin ks}{k}\right),
 \qquad s\in (-L_0(k)/2, L_0(k)/2).
\]
Let $A_0=A_0(k)$ be in Proposition \ref{eq:disk-volume-area} and consider $$X(A_0, s)=\tilde{X}(k, s)$$ as a variation of the model arcs, with variation parameter $A_0$.
Define $$u_k(s)=\left<\frac{\partial X}{\partial A_0}(s), N(s)\right>,$$ where $N$ is outward unit normal to $\gamma_k$.
\begin{prop}\label{lem:disk-model-speed}

\[
 u_k(s)=
 \frac{1-\cos(ks)/\sqrt{1+k^2}}{L_0(k)-2/(1+k^2)}.
\]
At $k=0$, it is \[u_0(s)=3(1+s^2)/8, \quad s\in [-1,1].\] In particular, $u_k\ge 0$ and 
\[
 \int_{\gamma_k}u_k\,\dd s=1.
\]
\end{prop}

\begin{proof}
The outward unit normal is $N=(-\cos ks,\sin ks)$. By direct computation, we see
Thus
\[
\left<\frac{\partial \tilde{X}}{\partial k}(s), N(s)\right>
 =-\frac{1-\cos(ks)/\sqrt{1+k^2}}{k^2}.
\]
Using
\[
 A_0'(k)=-\frac{L_0(k)-2/(1+k^2)}{k^2},
\]
this gives
\[
 u_k(s)=\frac{1}{A_0'(k)}\left<\frac{\partial \tilde{X}}{\partial k}(s), N(s)\right>
 =\frac{1-\cos(ks)/\sqrt{1+k^2}}{L_0(k)-2/(1+k^2)}.
\]
At $k=0$, Taylor expansion gives
\[
 \begin{aligned}
 1-\frac{\cos(ks)}{\sqrt{1+k^2}}
 &=\frac{1+s^2}{2}k^2+O(k^4),\\
 L_0(k)-\frac{2}{1+k^2}
 &=\frac43k^2+O(k^4).
 \end{aligned}
\]
Hence $u_0(s)=3(1+s^2)/8$ on $[-1,1]$.
From the first variation formula of area, we see that
\[
1= \frac{\partial A_0}{\partial A_0}=\int_{\gamma_k}\left<\frac{\partial X}{\partial A_0}, N\right>\,\dd s=\int_{\gamma_k} u_k \dd s.
\]
\end{proof}


\section{Second variation estimates}

The isoperimetric profile expresses the least boundary length as a
function of enclosed area. Under a variation of an isoperimetric
region, both area and length change. To obtain a sharp differential
inequality for the profile, we need to choose a suitable normal speed
$u$, so that the resulting estimate for $Q(u)$ allows a comparison with
the model profile.
Proposition~\ref{lem:disk-model-speed} shows that the normal speeds along the
model family of isoperimetric regions are not constant. We use these
model speeds as a guide to construct test functions on general
isoperimetric boundarys through the arclength parameter.

For $k\ge 0$ and $0<L\le L_0(k)$, where $L_0(k)$ is given in Proposition \ref{eq:disk-volume-area}. Define
\begin{equation}\label{eq:hatQ}
 \begin{aligned}
 \widehat Q(k,L)
 &=-\frac{k^2}{\displaystyle L-
 \frac{2\sin(kL/2)}{k\bigl(\cos(kL/2)+k\sin(kL/2)\bigr)}}
 &&(k>0),\\
 \widehat Q(0,L)&=\lim_{k\to 0+}\widehat Q(k,L)= -\frac{2}{L^2(1-L/6)}.
 \end{aligned}
\end{equation}
Define $u_{k, L}: [-L/2,L/2]\to \mathbb{R}$ by 
\begin{equation}\label{eq:arc-speed}
 \begin{aligned}
 u_{k, L}(s)&=
 \frac{\displaystyle 1-
 \frac{\cos(ks)}{\cos(kL/2)+k\sin(kL/2)}}
 {\displaystyle L-
 \frac{2\sin(kL/2)}{k\bigl(\cos(kL/2)+k\sin(kL/2)\bigr)}}
 &&(k>0),\\
 u_{0, L}(s)&=\lim_{k\to 0+}u_{k, L}(s)=\frac{s^2+L-L^2/4}{L^2(1-L/6)}.
 \end{aligned}
\end{equation}
$u_{k,L}$ is the solution to the following ODE:
\[
 -u''-k^2u=\widehat Q(k,L),\qquad
 u'(L/2)=u(L/2),\qquad
 u'(-L/2)=-u(-L/2).
\]
By direct computation and integration by parts, we see that $u=u_{k, L}$ satisfies
\[
 \int_{-L/2}^{L/2}u\,\dd s=1,
\]
and
\[
 \int_{-L/2}^{L/2}\bigl((u')^2-k^2u^2\bigr)\,\dd s
 -u(-L/2)^2-u(L/2)^2=\widehat Q(k,L)<0.
\]

\begin{rem}In view of Propositions \ref{eq:disk-volume-area} and \ref{lem:disk-model-speed}, we see the following:
    When $L=L_0(k)$, $$\pi\widehat Q(k,L_0(k))=h_{\mathbb{B}^2}''(t), \quad \hbox{ for }t=\frac{A_0(k)}{|\mathbb{B}^2|}$$  and
    $$u_{k,L_0(k)}=u_k,$$ where $u_k$ is given in Proposition \ref{lem:disk-model-speed}.
    Indeed, the definition of $\widehat Q(k,L)$ is motivated by that of $h_{\mathbb{B}^2}''(t)$.
\end{rem}

\begin{prop}\label{prop:arc-test}
$\widehat Q(k,L)$ is well defined and $\widehat Q(k,L)<0$. For fixed $k\ge0$, the function $L\mapsto L\widehat Q(k,L)$ is strictly
increasing on $(0,L_0(k)]$.
\end{prop}

\begin{proof}
First consider $k>0$. Set
\[
 x=\frac{kL}{2},\qquad f(x)=x(\cot x+k).
\]
Since $L\le L_0(k)$, we have $0<x\le\arctan k<\pi/2$, and hence
$k\ge\tan x$. Thus
\[
 f'(x)=\cot x+k-x\csc^2x
 \ge\frac{\tan x-x}{\sin^2x}>0.
\]
It follows that $f(x)>f(0+)=1$. this gives
\[
 L-\frac{2\sin(kL/2)}
 {k\bigl(\cos(kL/2)+k\sin(kL/2)\bigr)}
 =L\left(1-\frac1{f(x)}\right)>0.
\]
Therefore $u_{k,L}$ and $\widehat Q(k,L)$ are well defined, and
$\widehat Q(k,L)<0$.

Note
\[
 L\widehat Q(k,L)=-k^2\frac{f(x)}{f(x)-1}.
\]
Since $\dd x/\dd L=k/2$, we have
\[
 \frac{\dd}{\dd L}\bigl(L\widehat Q(k,L)\bigr)
 =\frac{k^3f'(x)}{2(f(x)-1)^2}>0.
\]

Next consider $k=0$. We have $0<L\le L_0(0)=2$, so $\widehat Q(0,L)<0$. By direct computation
\[
 \frac{\dd}{\dd L}\bigl(L\widehat Q(0,L)\bigr)
 =\frac{2(1-L/3)}{(L-L^2/6)^2}>0.
\]
\end{proof}

\begin{prop}Let $(M^2, g)$ be a compact surface with boundary such that $K_M\ge 0$ and $\kappa_{\partial M}\ge 1$.
Let $\gamma$ be a free-boundary arc in $M$ with constant  geodesic curvature
$k$, and denote $|\gamma|=L$.  Let $u_{k, L}$ be defined as \eqref{eq:arc-speed}. Then
 $$Q(u_{k,L})\le\widehat Q(k,L).$$
\end{prop}

\begin{proof}
For simplicity, we write $u=u_{k, L}$. Since $K_M\ge 0$ and $\kappa_{\partial M}\ge 1$, we have
\[
 \begin{aligned}
 Q(u)
 &=\int_\gamma\bigl((u')^2-(k^2+K_M)u^2\bigr)\,\dd s
 -\sum_{p\in\partial\gamma}\kappa_{\partial M}(p)u(p)^2\\
 &\le\int_{-L/2}^{L/2}\bigl((u')^2-k^2u^2\bigr)\,\dd s
 -u(-L/2)^2-u(L/2)^2\\
 &=\widehat Q(k,L).
 \end{aligned}
\]
\end{proof}

\section{Proof of the main theorem}

\begin{proof}[Proof of the {\hyperref[thm:main]{Main Theorem}}]
Both profiles $h_M$ and $h_{\mathbb{B}^2} $are
continuous on $[0,1]$ and vanish at the endpoints.
If  $h_M-h_{\mathbb{B}^2}$ attains its minimum at $0$ or $1$, then $h_M-h_{\mathbb{B}^2}\ge 0$.

Next we consider the case $h_M-h_{\mathbb{B}^2}$ attains its minimum at some $t_0\in (0, 1)$, then $h_M(t_0)\le h_{\mathbb{B}^2}(t_0)$.
By symmetry,
we may assume that $0<t_0\le1/2$. Let $k_0\ge 0$ be such that $A_0(k_0)=\pi t_0$.

Let $E$ be an isoperimetric region with $|E|=t_0|M|$, and 
$\gamma$ be its relative boundary, with geodesic curvature $\kappa_\gamma$.
Take a smooth upper support $\widehat h_M$ given in
Proposition~\ref{lem:upper-support}. Near $t_0$,
\[
 \begin{aligned}
 \widehat h_M(t)-h_{\mathbb{B}^2}(t)
 &\ge h_M(t)-h_{\mathbb{B}^2}(t)
 \ge h_M(t_0)-h_{\mathbb{B}^2}(t_0)
 =\widehat h_M(t_0)-h_{\mathbb{B}^2}(t_0).
 \end{aligned}
\]
Thus $\widehat h_M-h_{\mathbb{B}^2}$ has a local minimum at $t_0$, and
\[
 0=(\widehat h_M-h_{\mathbb{B}^2})'(t_0)
 =\kappa_\gamma-k_0.
\]
In particular,
\[
 \kappa_\gamma=k_0\ge0.
\]

By Proposition \ref{prop:regularity}, $\gamma$ is connected.
We first show that it meets $\partial M$. Suppose instead that
$\partial\gamma=\varnothing$. Then $\gamma$ is a simple closed
curve. Proposition~\ref{lem:geometry-two} excludes $\kappa_\gamma=0$, so $\kappa_\gamma>0$.

Since $M$ is a disk, $\gamma$ bounds an inner disk
$U\subset\operatorname{int}M$. 
Up to a set of measure zero, $E$ is either $U$ or $M\setminus\overline U$. If $E=M\setminus\overline U$,
the curvature of $\gamma$ with respect to the outward normal of
$U$ is $-\kappa_\gamma$. Using the Gauss--Bonnet theorem applying on $U$ and
Proposition~\ref{lem:geometry-two}, we see
\[
 -\kappa_\gamma|\gamma|=2\pi-\int_{U} K_M\,\dd \mu>0,
\]
contradicting $\kappa_\gamma>0$. Hence $E=U$.

 Applying the Gauss--Bonnet theorem to $U$ and $M$, using $K_M\ge 0$ and $\kappa_{\partial M}\ge 1$, and Proposition~\ref{lem:geometry-two}, we obtain
\[
 \begin{aligned}
 \kappa_\gamma|\gamma|
 &=2\pi-\int_U K_M\,\dd \mu=\int_{M\setminus\overline U} K_M\,\dd \mu+\int_{\partial M}\kappa_{\partial M}\,\dd s\ge|\partial M|\ge2|M|.
 \end{aligned}
\]
 Since $|\gamma|=|M|h_M(t_0)$,
this implies
\[
 \kappa_\gamma h_M(t_0)\ge2.
\]
On the other hand, $h_M(t_0)\le h_{\mathbb{B}^2}(t_0)$, and
Proposition~\ref{prop:disk-model} gives
\[
 \kappa_\gamma h_M(t_0)
 \le\kappa_\gamma h_{\mathbb{B}^2}(t_0)=k_0h_{\mathbb{B}^2}(t_0)
 =\frac{2\arctan k_0}\pi<1.
\]
This contradiction shows that $\partial\gamma\ne\varnothing$.
Hence $\gamma$ is a free-boundary curve.


Write $L=|\gamma|$. By Proposition~\ref{lem:geometry-two} and the choice
of $t_0$,
\[
 0<L=|M|h_M(t_0)
 \le \pi h_{\mathbb{B}^2}(t_0)=L_0(k_0).
\]
Since $\kappa_\gamma=k_0$,
We may therefore take the normalized speed $u_{k_0, L}$ from
Proposition~\ref{prop:arc-test}, which satisfies
\[
 \int_\gamma u_{k_0, L}\,\dd s=1,
 \qquad Q(u_{k_0, L})\le\widehat Q(k_0,L).
\]

Let $\widehat h_M$ be the upper support corresponding to $u_{k_0, L}$.
As before, $\widehat h_M-h_{\mathbb{B}^2}$ has a local minimum at $t_0$.
Hence
\begin{eqnarray}\label{xeq1}
 0\le(\widehat h_M-h_{\mathbb{B}^2})''(t_0)
 =|M|Q(u_{k_0, L})-h_{\mathbb{B}^2}''(t_0).
\end{eqnarray}
On the other hand, since $L\mapsto L\widehat Q(k_0,L)$ is strictly
increasing by Proposition~\ref{prop:arc-test}, we obtain
\begin{equation}\label{xeq2}
 |M|Q(u_{k_0, L})
 \le\frac{L\widehat Q(k_0,L)}{h_M(t_0)}\le\frac{L_0(k_0)\widehat Q(k_0,L_0(k_0))}{h_M(t_0)}
 =\frac{h_{\mathbb{B}^2}(t_0)h_{\mathbb{B}^2}''(t_0)}{h_M(t_0)}
 \le h_{\mathbb{B}^2}''(t_0),
 \end{equation}
where the last inequality uses $h_{\mathbb{B}^2}(t_0)\ge h_M(t_0)$ and
$h_{\mathbb{B}^2}''(t_0)<0$. By virtue of \eqref{xeq1} and \eqref{xeq2}, we see that all inequalities in \eqref{xeq2} are equalities.
In particular, $$h_M(t_0)=h_{\mathbb{B}^2}(t_0).$$Since $t_0$ is minimum point of $h_M-h_{\mathbb{B}^2}$, we conclude that
$h_M\ge h_{\mathbb{B}^2}$ for any $t\in [0, 1]$.
This finishes the proof of first assertion in Theorem \ref{thm:main}.

We now consider the rigidity. Suppose that
\[
 h_M(t_0)=h_{\mathbb{B}^2}(t_0)
\]
for some $t_0\in(0,1)$. The above argument is still valid.  From equality in \eqref{xeq2}, we also obtain 
$L=L_0(k_0)$, since
$L\mapsto L\widehat Q(k_0,L)$ is strictly increasing. Consequently
\[
 |M|=\frac{L}{h_{\mathbb{B}^2}(t_0)}=\frac{L_0(k_0)}{h_{\mathbb{B}^2}(t_0)}=\pi.
\]

By using the Heintze-Karcher inequality (see Ros \cite[Theorem 1]{Ros}) and Proposition \ref{lem:geometry-two}, we obtain
\[
 \pi=|M|
 \le\frac12\int_{\partial M}\frac1\kappa\,\dd s
 \le\frac{|\partial M|}{2}\le\pi.
\]
Thus equality holds in the Heintze-Karcher inequality.
Hence $M$ is isometric to a Euclidean unit disk. The proof of Theorem \ref{thm:main} is completed. 
\end{proof}

\end{document}